\documentstyle[12pt,amscd,amssymb]{amsart}
\newtheorem{thm}{Theorem}[section]

\newtheorem{prop}[thm]{Proposition}

\begin{document}

\newcommand{\C}{{\Bbb C}}
\newcommand{\HH}{{\Bbb H}}
\newcommand{\M}{{\cal M}}
\newcommand{\R}{{\Bbb R}}
\newcommand{\Z}{{\Bbb Z}}
\newcommand{\E}{{\cal E}}
\renewcommand{\O}{{\cal O}}
\renewcommand{\P}{{\Bbb P}}
\newcommand{\Q}{{\cal Q}}
\newcommand{\rk}{\o{rank}}
\title
{Moduli of Vector Bundles on Curves in Positive Characteristics}
\author{Kirti Joshi \and Eugene Z. Xia}
\date{\today} 
\address{ 
Department of Mathematics,
University of Arizona, Tucson, AZ 85721}
\email{kirti@@math.arizona.edu {\it (Joshi)}, exia@@math.arizona.edu {\it (Xia)}}
\maketitle
\begin{abstract}
Let $X$ be a projective curve of genus 2
over an algebraically closed field of characteristic 2.
The Frobenius map on $X$ induces a rational map
on the moduli scheme of rank-2 bundles.  We show that up
to isomorphism, there is only one (up to tensoring
by an order two line bundle) semi-stable vector 
bundle of rank 2 (with determinant equal to a theta characteristic)
whose Frobenius pull-back is not semi-stable.  The indeterminacy of
the Frobenius map at this point can be resolved by introducing
Higgs bundles.
\end{abstract}

\section{Introduction and Results}
Let $X$ be a smooth projective curve of genus 2 over an algebraically
closed field $k$ of characteristic $p > 0$.  
Let $\Omega$ be its canonical bundle.
Define the (absolute)
Frobenius morphism \cite{Ka1, Ka2}
\begin{equation*}
F : X \longrightarrow X
\end{equation*}
which maps local sections $f \in \O_X$ to $f^p$.
As $X$ is smooth, $F$ is a (finite) flat map.

Let $J^0, J^1$ be the moduli schemes of isomorphism 
classes of line bundles of degree $0$ and $1$, respectively.
Choose a theta characteristic $L_{\theta} \in J^1$.
Denote by $S_O$ (resp. $S_{\theta}$) the moduli scheme of
S-equivalence classes of semi-stable vector bundles of rank
$2$ and determinant $O_X$ (resp. $L_{\theta}$) on $X$ \cite{Se}.
We study the Frobenius pull-backs of the bundles in $S_O$ and 
$S_{\theta}$.
The geometry of $S_{\theta}$ has been studied extensively
by Bhosle \cite{Bh}.  

The operation of Frobenius pull-back has a tendency to
destabilize bundles \cite{Ra}.  In particular, the map 
$V \longmapsto F^*(V)$ is rational on the moduli scheme.

The Frobenius destabilizes only finite many bundles in $S_O$
(see Theorem~\ref{prop:3.2b}).  For any $V \in S_O$, 
Proposition~\ref{prop:3.5} gives a necessary and sufficient
criterion for $F^*(V)$ to be non-semi-stable in terms of theta
characteristic.

For a given vector bundle $V$ on $X$, let
$$
J_2(V) = \{V \otimes L : L \in J^0, L^2 = O_X\}
$$
\begin{thm} \label{thm:2}
Suppose $p=2$.
Then there exists a bundle $V_1 \in Ext^1(L_{\theta}, \O_X)$ such
that if
$$
V \in S_{\theta} \setminus J_2(V_1),
$$
then $F^*(V)$ is semi-stable.
Hence, the Frobenius map induces a map
$$
\Omega^{-1} \otimes F^* : S_{\theta} \setminus J_2(V_1) \longrightarrow S_O.
$$
\end{thm}

We show that there is a natural way of resolving the indeterminacy
of the Frobenius map at the points in $J_2(V_1)$, by replacing 
$S_O, S_{\theta}$ with moduli schemes of suitable Higgs bundles.
Denote by $S_O(\Omega)$ (resp. $S_{\theta}(L_{\theta})$) the moduli scheme
of semi-stable Higgs bundles with associated line bundle $O_X$ 
(resp. $L_{\theta}$) \cite{Ni}.  
For any Higgs bundle on $X$, one may also consider its
Frobenius pull-back.  
\begin{thm} \label{thm:3}
Suppose $p=2$.
\begin{enumerate}
\item If $(V, \phi) \in S_{\theta}(L_{\theta})$, then either $V \in S_{\theta}$
or $V \in J_2(O_X \oplus L_{\theta})$.
\item There exist Higgs fields $\phi_0$ and $\phi_1$
such that
$(F^*(W), F^*(\phi_0))$ and $(F^*(V), F^*(\phi_1))$
are semi-stable for all
$W \in J_2(O_X \oplus L_{\theta})$ and $V \in J_2(V_1)$.
\end{enumerate}
Hence, the Frobenius defines a map on
a Zariski open set $U \subset S_{\theta}(L_{\theta})$ 
$$
\Omega^{-1} \otimes F^* : U \longrightarrow S_O(\O_X),
$$
where $U$ contains the scheme $S_{\theta} \setminus J_2(V_1)$ and the points
$(W, \phi_0), (V, \phi_1)$ for any $W \in J_2(O_X \oplus L_{\theta})$ and 
$V \in J_2(V_1)$.
\end{thm}

Cartier's theorem gives a criterion for descent under
Frobenius \cite{Ka1}.  Higgs bundles appear
naturally in characteristic $p > 0$ context.  
To see this, let $(V, \nabla)$
be a vector bundle with a (flat) connection,
$$
\nabla : V \longrightarrow \Omega \otimes_{O_X} V.
$$
One associates to the pair $(V, \nabla)$ its $p$-curvature 
which is a homomorphism of $O_X$-modules \cite{Ka1, Ka2}:
$$
\psi : V \longrightarrow F^*(\Omega) \otimes_{O_X} V.
$$
Thus the pair $(V, \nabla)$ gives a Higgs bundle with associated line bundle
$F^*(\Omega)$.  

\centerline{\sc Acknowledgments}
We thank Professor Usha Bhosle for reading a previous 
version and for her comments and suggestions for improvement. 
We thank Professors Minhyong Kim, N. Mohan Kumar, V. B. Mehta and S. Ramanan
for insightful discussions and comments.  Finally, we thank the referee
for his or her comments.

\section{Bundle Extensions and the Frobenius Morphism}
Suppose $L$ is a line bundle on $X$.  Then $F^*(L) = L^p$.
The push-forward, $F_*(O_X)$, is a vector bundle of rank $p$
and one has the exact sequence of vector bundles \cite{Ra}
\begin{equation*}
0 \longrightarrow O_X \longrightarrow F_*(O_X)
\longrightarrow B_1 \longrightarrow 0.
\end{equation*}
Tensoring the sequence with a line bundle $L$ and using the projection
formula, we obtain
\begin{equation*}
0 \longrightarrow L \longrightarrow F_*(L^p)
\longrightarrow B_1 \otimes L \longrightarrow 0.
\end{equation*}
The associated long cohomology sequence is
\begin{equation*}
\cdots \longrightarrow H^0(B_1 \otimes L) \longrightarrow H^1(L)
\stackrel{f_L}{\longrightarrow} H^1(F_*(L^p)) \longrightarrow \cdots
\end{equation*}
Since $F$ is an affine morphism, the Leray spectral sequence
for $F$ degenerates at $E_2$.
Hence 
\begin{equation*}
H^i(F_*(L^p)) \cong H^i(L^p).
\end{equation*}
Substituting this into the long exact sequence, one obtains
\begin{equation}
\cdots \longrightarrow H^0(B_1 \otimes L) \longrightarrow H^1(L)
\stackrel{f_L}{\longrightarrow} H^1(L^p) \longrightarrow \cdots
\end{equation}
Suppose $V \in Ext^1(L_2, L_1) \cong H^1(L_2^{-1} \otimes L_1)$, i.e.
\begin{equation*} \label{eqn:1}
0 \longrightarrow L_1 \longrightarrow V
\longrightarrow L_2 \longrightarrow 0,
\end{equation*}
where $L_1, L_2$ are line bundles.
Since $F$ is a flat morphism, we have
\begin{equation*}
0 \longrightarrow F^*(L_1) \longrightarrow F^*(V)
\longrightarrow F^*(L_2) \longrightarrow 0.
\end{equation*}
This gives a map
\begin{equation*}
F^* : Ext^1(L_2, L_1) \longrightarrow Ext^1(F^*(L_2), F^*(L_1)) \cong
Ext^1(L_2^p, L_1^p).
\end{equation*}
Take $L = L_2^{-1} \otimes L_1$ in (1).
\begin{prop} \label{prop:2.1}
$F^*(V) = L_1^p \oplus L_2^p$ if and only if $V$ is in the
image of the connecting homomorphism 
$$
H^0(B_1 \otimes L_2^{-1} \otimes L_1)
\longrightarrow H^1(L_2^{-1} \otimes L_1).
$$
\end{prop}
\begin{pf}
Since the functors $\Gamma(X,.)$
and $Hom(O_X,.)$ are equivalent,
the diagram
\begin{equation*}
\begin{CD}
H^1(L_2^{-1} \otimes L_1)  @>f_{L_2^{-1} \otimes L_1}>>       H^1(L_2^{-p} \otimes L_1^p)\\
@VV{id}V                                 @VV{id}V\\
Ext^1(L_2, L_1)            @>F^*>>     Ext^1(L_2^p, L_1^p)
\end{CD}
\end{equation*}
commutes.  Now the proposition follows directly from the
long exact sequence
$$
\cdots \longrightarrow H^0(B_1 \otimes L_2^{-1} \otimes L_1)
\longrightarrow H^1(L_2^{-1} \otimes L_1)
\longrightarrow H^1(L_2^{-p} \otimes L_1^p) 
\longrightarrow \cdots
$$
\end{pf}

\section{The moduli of Semi-Stable Vector and Higgs Bundles}
Suppose $V$ is a vector bundle on $X$.  The slope of
$V$ is defined as
$$
\mu(V) = \deg(V) / \mbox{rank}(V).
$$
A vector bundle $V$ is semi-stable (resp. stable)
if for every proper subbundle $W$ of $V$, $\mu(W) \le \mu(V)$ 
(resp. $\mu(W) < \mu(V)$).  The schemes $S_O$ and $S_{\theta}$ are defined to
be the moduli schemes of all $S$-equivalence classes \cite{Se} of rank 2 
semi-stable vector bundles with determinant equal to 
$O_X$ and $L_{\theta}$, respectively.

A Higgs bundle $(V,\phi)$ with an associated line bundle $L$ 
on $X$ consists of a vector bundle $V$
and a Higgs field which is a morphism of bundles:
$$
\phi : V \longrightarrow V \otimes L.
$$
Frobenius pulls back Higgs fields
$$
F^*(V) \stackrel{F^*(\phi)}{\longmapsto} F^*(V) \otimes F^*(L),
$$
hence, pulls back Higgs bundles.

A Higgs bundle $(V,\phi)$ is said to be semi-stable (resp. stable)
if for every proper subbundle $W$ of $V$, satisfying  
$\phi(W) \subset W \otimes L$, one has 
$\mu(W) \le \mu(V)$ (resp. $\mu(W) < \mu(V)$).
The scheme $S_O(\Omega)$ (resp. $S_{\theta}(L_{\theta})$) is defined to
be the moduli scheme of all $S$-equivalence 
classes of rank 2 semi-stable Higgs 
bundles on $X$ with determinant $O_X$ (resp. $L_{\theta}$) and with
associated line bundle 
$\Omega$ (resp. $L_{\theta}$) \cite{Ni}.

Let
$$
K = \{V \in S_O : V \mbox{ is semi-stable but not stable}\}.
$$
Suppose $V \in K$.
Then there exists $L \in J^0$ such that
$$
0 \longrightarrow L^{-1} \longrightarrow V \longrightarrow L \longrightarrow 0.
$$
The pull-back of $V$ by Frobenius then fits into the following
sequence
$$
0 \longrightarrow L^{-p} \stackrel{f_1}{\longrightarrow} F^*(V)
\stackrel{f_2}{\longrightarrow} L^p \longrightarrow 0.
$$
\begin{prop} \label{prop:3.1}
$$
F^* : K \longrightarrow K.
$$
is a well-defined morphism.
\end{prop}
\begin{pf}
Let $H \subset V$ be a subbundle of maximum degree.
If $f_2 |_H = 0$, then $H = L^{-p}$ and
$\deg(H) = \deg(L^{-p}) = 0$.  If $f_2 |_H \neq 0$,
then $\deg(H) \le \deg(L^p) = 0$.
\end{pf}

In general, $F^*(V)$ may not be semi-stable.
For example,
a theorem of Raynaud states that the bundle $B_1$ is always semi-stable while
$F^*(B_1)$ is never semi-stable for all $p > 2$ \cite{Ra}.  
The following theorem 
was communicated to Joshi by V.B. Mehta:
\begin{thm}\label{prop:3.2b}
Let $X$ be a curve of genus $2$ over an algebraically closed
field of characteristic $p > 2$.
Then there exists a finite set S,
such that $F^*(V)$ is semi-stable for all $V \in S_O \setminus S$.
In other words, $F^*$ induces a morphism:
$$
F^* : S_O \setminus S \longrightarrow S_O.
$$
\end{thm}
\begin{pf}
By a theorem of Narasimhan-Ramanan, when $p=0$ \cite{Na},
$S_O \cong {\mathbb P}^3$.  Moreover, as was remarked to one of
us by Ramanan, the proof given there works
in all characteristic $p \neq 2$.
The Frobenius morphism is defined on a non-empty Zariski open set $U$ in
$S_O \cong {\mathbb P}^3$.
By Proposition~\ref{prop:3.1}, $U$ contains $K$ which is an ample divisor 
in ${\mathbb P}^3$.
Therefore $S_O \setminus U$ is of co-dimension 3, 
hence, is a finite set.  Note that $K$ can also be
identified with the Kummer surface of $J^0$ in ${\mathbb P}^3$ \cite{Na}.
\end{pf}

When $X$ is ordinary, $F^*$ is \'{e}tale on a non-empty
Zariski open set of $S_O$ \cite{Me}.
Although unable to identify explicitly this finite set
upon which the Frobenius is not defined, we provide the 
following criterion.
\begin{prop} \label{prop:3.5}
Let $X$ be a curve of genus $2$ over an algebraically closed
field of characteristic $p > 0$.
Suppose $V \in S_O$ Then $F^*(V)$ is
not semi-stable if and only if $F^*(V)$ is an extension
$$
0 \longrightarrow M \longrightarrow F^*(V) \longrightarrow M^{-1} 
\longrightarrow 0,
$$
where $M \in J_2(L_{\theta})$.
\end{prop}

\begin{pf}
One direction is clear.
We use inseparable
descent to prove the other direction. Suppose $F^*(V)$ is not
semi-stable.  Then we have an exact sequence
$$ 
0 \longrightarrow M \longrightarrow F^*(V) \longrightarrow M^{-1} \to 0,
$$
where $\deg(M) > 0$. 

Following \cite{Ka1}, consider the natural connection on $F^*(V)$
with zero $p$-curvature. Then the second fundamental form of this
connection is a morphism
\begin{equation*}
T_X \longrightarrow Hom(M,M^{-1})=M^{-2}.
\end{equation*}
As $V$ is semi-stable, this
morphism must not be the zero morphism.  In other words,
$M^{-2} \otimes \Omega$ has a non-zero section.  Since $\deg(M) > 0$
and $\deg(\Omega) = 2$, we must have 
$\Omega = M^2$.  Hence $M \in J_2(L_{\theta})$.
\end{pf}

\section{The Moduli Spaces in Characteristic 2}
In this section, we assume $p=2$.
Then $B_1$ is a line bundle and equal to 
a theta characteristic \cite{Ra}.  Choose $L_{\theta}$ 
to be $B_1$.

\subsection{The moduli of semi-stable bundles}
Suppose $V \in S_{\theta}$. 
By a theorem in \cite{Na},
there exist
$L_1 \in J^0, L_2 \in J^1$ with $L_1 \otimes L_2 = L_{\theta}$ such that
$V \in Ext^1(L_2, L_1)$.
Since $L_{\theta} = B_1$, $h^0(B_1 \otimes L_2^{-1} \otimes L_1)$ is
1 if $L_{\theta} = L_2 \otimes L_1^{-1}$ and $0$ otherwise.  Hence,
by Proposition~\ref{prop:2.1}, 
there is a unique (up to a scalar) $V_1$ not isomorphic to 
$O_X \oplus L_{\theta}$ and
$$
0 \longrightarrow O_X \longrightarrow V_1 \longrightarrow L_{\theta} 
\longrightarrow 0
$$ 
such that $F^*(V_1) = O_X \oplus \Omega$.  
It is immediate that $V_1$ is stable \cite{Na}.

Suppose $V \not\in J_2(V_1)$.
Then by Proposition~\ref{prop:2.1},
$$
F^*(V) \neq F^*(L_1) \oplus F^*(L_2) = L_1^2 \oplus L_2^2.
$$
If $M \subset F^*(V)$ is a destabilizing subbundle, i.e. $\deg(M) \ge 2$, 
then $M^{-1} \otimes L_2^2$ has a global section implying 
$$
\deg(M) \le \deg(L_2^2) = 2.
$$
Moreover if $\deg(M) = 2$,
then $M = L_2^2$ implying that $F^*(V)$ contains 
$L_2^2$ as a subbundle.
Then the sequence
$$
0 \longrightarrow L_1^2 \longrightarrow F^*(V)
\longrightarrow L_2^2 \longrightarrow 0,
$$
splits.  This is a contradiction.
This proves Theorem~\ref{thm:2}.

\subsection{Restoring Frobenius Stability: Higgs Bundles}
The scheme $S_{\theta}$ embeds in $S_{\theta}(L_{\theta})$
by the map $V \longmapsto (V,0)$.  
If $(V, \phi) \in S_{\theta}(L_{\theta})$
and $V$ is not semi-stable, then $V$ is an extension
\begin{equation}
0 \longrightarrow L_1 \longrightarrow V 
\stackrel{f}{\longrightarrow} L_2 \longrightarrow 0,
\end{equation}
where $\deg(L_1) \ge 1 > \deg(L_2)$.  Moreover $\phi(L_1)$
is not contained 
in $L_1 \otimes L_{\theta}$
(otherwise $\phi(L_1) \subset L_1 \otimes L_{\theta}$
implying $(V,\phi)$ is not semi-stable).
This implies that
there exists a line bundle $H \subset V$ such that $H \neq L_1$ and 
$\phi(L_1) \subset H \otimes L_{\theta}$.  Then
$$
\deg(L_1) \le \deg(H) + \deg(L_{\theta}).
$$
Since $L_1 \neq H$, $0 \neq f(H) \subset L_2$ implies that 
$\deg(H) \le \deg(L_2)$.  To summarize, we have the following inequalities:
$$
\deg(L_2) + \deg(L_{\theta}) \ge \deg(H) + \deg(L_{\theta}) \ge \deg(L_1) > \deg(L_2).
$$
Since $\deg(L_{\theta}) = 1$, $\deg(L_1) = \deg(H) + 1 = \deg(L_2) + 1 = 1$.  
The degree of $H$ is thus zero implying that $f(H) = L_2$, 
so the exact sequence (2) splits.  
In addition, since $0 \neq \phi(L_1) \subset H \otimes L_{\theta}$, 
$\phi |_{L_1}$
must be a non-zero constant morphism and
$$
L_1 = L_2 \otimes L_{\theta}.
$$
Since $L_1 \otimes L_2 = L_{\theta}$,
$V \in J_2(O_X \oplus L_{\theta})$.  This proves the first part of 
Theorem~\ref{thm:3}.

Suppose $(V,\phi) \in S_{\theta}(L_{\theta})$.
If $V \in S_{\theta} \setminus J_2(V_1)$, then $F^*(V)$ is semi-stable
by Theorem~\ref{thm:2}; hence,
$(F^*(V), F^*(\phi))$ is semi-stable.

\noindent{\sl The split case:} Suppose $W = L \oplus (L \otimes L_{\theta})$,
where $L \in J_2(O_X)$.

We take
the Higgs field $\phi_0$ to be the identity map:
$$
1 = \phi_0 : L \otimes L_{\theta} \longrightarrow L \otimes L_{\theta}.
$$
If $M \subset W$,
then either $M = L \otimes L_{\theta}$ or $\mu(M) < \mu(W)$.
Since $L \otimes L_{\theta}$ is not $\phi_0$-invariant, 
$(W,\phi_0)$ is stable.
The Frobenius pull-back  $F^*(\phi_0)$ is again a constant map
$$
F^*(\phi_0) : \Omega \longrightarrow O_X \otimes \Omega.
$$
Now if $N \subset O_X \oplus \Omega$, then either 
$N = \Omega$ or $\mu(N) < \mu(O_X \oplus \Omega)$.
Since $\Omega$ is not $F^*(\phi_0)$-invariant, 
$(F^*(W), F^*(\phi_0))$ is stable.

\noindent{\sl The non-split case:} Suppose $V = L \otimes V_1$, where
$L \in J_2(O_X)$.

The bundle $V$ is a non-trivial extension:
\begin{equation}
0 \longrightarrow L \stackrel{f_1}{\longrightarrow} V 
\stackrel{f_2}{\longrightarrow} L \otimes L_{\theta} \longrightarrow 0.
\end{equation}
Tensoring the sequence with $L_{\theta}$ gives
\begin{equation}
0 \longrightarrow L \otimes L_{\theta} \stackrel{g_1}{\longrightarrow} V
\otimes L_{\theta} 
\stackrel{g_2}{\longrightarrow} L \otimes \Omega \longrightarrow 0.
\end{equation}
Set 
$$
\phi_1 = g_1 \circ \phi_0 \circ f_2 : V \longrightarrow L_{\theta} \otimes V.
$$
The Frobenius pull-back decomposes $V$:
$$
F^*(V) = O_X \oplus \Omega.
$$
Pulling back the exact sequences (3) and (4) by Frobenius gives
$$
0 \longrightarrow O_X \stackrel{F^*(f_1)}{\longrightarrow} O_X \oplus \Omega 
\stackrel{F^*(f_2)}{\longrightarrow} \Omega \longrightarrow 0
$$
$$
0 \longrightarrow O_X \otimes \Omega \stackrel{F^*(g_1)}{\longrightarrow} 
(O_X \oplus \Omega) \otimes \Omega
\stackrel{F^*(g_2)}{\longrightarrow} \Omega \otimes \Omega \longrightarrow 0
$$
Suppose $N \subset O_X \oplus \Omega$.  Then either $N = \Omega$ or
$\mu(N) < \mu(O_X \oplus \Omega)$.  
The Frobenius pull-back of $\phi_1$ is a composition:
$$
F^*(\phi_1) = F^*(g_1) \circ F^*(\phi_0) \circ F^*(f_2).
$$
Since the map $F^*(f_2)$ is surjective, the restriction map
$F^*(f_2)|_{\Omega}$ is an isomorphism.
The map $\phi_0$ is an isomorphism and $g_1$ is injective; hence,
$g_1 \circ \phi_0$ is injective.  This implies $F^*(g_1) \circ F^*(\phi_0)$
is injective.  Therefore $F^*(\phi_1)|_{\Omega}$ is injective.  
Since $\deg(\Omega) < \deg(\Omega \otimes \Omega)$, $F^*(\phi_1)|_{\Omega}$
being injective implies
$$
F^*(\phi_1)(\Omega) \not\subset \Omega \otimes \Omega \subset 
(O_X \oplus \Omega) \otimes \Omega.
$$
In other words,
$\Omega \subset O_X \oplus \Omega$ is not $F^*(\phi_1)$-invariant.
Hence $(F^*(V), F^*(\phi_1))$ is stable.  This proves Theorem~\ref{thm:3}.

\end{document}